\documentclass[letterpaper, 10 pt, conference]{ieeeconf}  

\title{\LARGE \bf
  Wasserstein Contraction Bounds on Closed Convex Domains with Applications to Stochastic Adaptive Control
}

\author{Tyler Lekang and Andrew Lamperski
}

\usepackage[style=ieee,backend=bibtex]{biblatex}
\usepackage{amsmath}
\usepackage{amssymb}
\usepackage[scr=esstix]{mathalpha}
\usepackage{bm}
\usepackage{bbm} 
\usepackage{hyperref}
\hypersetup{colorlinks,
	linkcolor=blue,
	citecolor=blue,
	urlcolor=magenta,
	linktocpage,
	plainpages=false}
\usepackage{cancel}
\usepackage{graphicx}
\usepackage{mathtools}
\newtheorem{theorem}{Theorem}
\newtheorem{corollary}{Corollary}

\newtheorem{remark}{Remark}

\newcommand{\E}{\mathbb{E}}
\renewcommand{\P}{\mathbb{P}}
\newcommand{\Q}{\mathbb{Q}}
\newcommand{\R}{\mathbb{R}}
\newcommand{\Rnm}{\R^{n+m}}

\newcommand{\ub}{\bm{u}}
\newcommand{\vb}{\bm{v}}
\newcommand{\Vb}{\bm{V}}
\newcommand{\wb}{\bm{w}}
\newcommand{\rb}{\bm{r}}
\newcommand{\xb}{\bm{x}} 
\newcommand{\yb}{\bm{y}}
\newcommand{\zb}{\bm{z}}
\newcommand{\mb}{\bm{m}}
\newcommand{\mub}{\bm{\mu}}

\newcommand{\psib}{\bm{\psi}}
\newcommand{\taub}{{\bm{\tau}}}
\newcommand{\Kb}{\bm{K}}
\newcommand{\thetab}{{\bm{\theta}}}
\newcommand{\Thetab}{{\bm{\Theta}}}
\newcommand{\Omegab}{\bm{\Omega}}
\newcommand{\Lb}{\bm{L}}

\newcommand{\thetah}{{\hat{\theta}}}

\newcommand{\Thetabar}{{\bar{\Theta}}}
\newcommand{\thetahb}{{\bm{\hat{\theta}}}}

\newcommand{\Ac}{\mathcal{A}}
\newcommand{\Lc}{\mathcal{L}}
\newcommand{\Kc}{\mathcal{K}}
\newcommand{\Xc}{\mathcal{X}}

\newcommand{\Hs}{\mathscr{H}}
\newcommand{\Cs}{\mathscr{C}}

\newcommand{\Rs}{\mathscr{R}}

\newcommand{\Vs}{\mathscr{V}}

\newcommand{\xw}{\widetilde{x}}

\newcommand{\xbw}{\widetilde{\xb}}
\newcommand{\zbw}{\widetilde{\zb}}

\newcommand{\Thetaw}{\widetilde{\Theta}}

\renewcommand{\d}{\text{d}} 
\newcommand{\dt}{\text{d}t} 

\newcommand{\tr}{\text{tr}}
\newcommand{\TV}{\mathrm{TV}}

\newcommand{\Ind}{\mathbb{I}}

\newcommand{\sigx}{G_x}

\begin{document}

\maketitle
\thispagestyle{empty}
\pagestyle{empty}

\begin{abstract}
This paper is motivated by the problem of quantitatively bounding the convergence of adaptive control  methods for stochastic systems to a stationary distribution. Such bounds are useful for analyzing statistics of trajectories and determining appropriate step sizes for simulations. To this end, we extend a methodology from (unconstrained) stochastic differential equations (SDEs) which provides contractions in a specially chosen Wasserstein distance. This theory focuses on unconstrained SDEs with fairly restrictive assumptions on the drift terms. Typical adaptive control schemes place constraints on the learned parameters and their update rules violate the drift conditions. To this end, we extend the contraction theory to the case of constrained systems represented by reflected stochastic differential equations and generalize the allowable drifts. We show how the general theory can be used to derive quantitative contraction bounds on a nonlinear stochastic adaptive regulation problem.
\end{abstract}

\section{Introduction}
Adaptive control has a rich history in the controls literature \cite{lavretsky2013robust},
\cite{hovakimyan2010L1}, \cite{ioannou2006adaptive}, and \cite{khalil2017high}. It has wide applications in areas such as robotics \cite{lewis1998neural}, aerospace systems \cite{lavretsky2013robust}, and electromechanical systems \cite{astolfi2007nonlinear}.

The typical approach utilizes Lyapunov-based design to update the parameters while guaranteeing stability.

In recent years, there has been a drive to connect adaptive control methods with techniques from reinforcement learning \cite{vrabie2013optimal,kamalapurkar2018reinforcement,westenbroek2020adaptive,boffi2020regret}. In parallel, methods from reinforcement learning have seen an explosion of work on linear giving precise optimality guaranatees \cite{recht2019tour,simchowitz2020improper,hazan2020nonstochastic}. These works rely on precise convergence bounds that are fairly straightforward for linear systems, but substantially more complex in stochastic nonlinear systems.

Convergence of stochastic nonlinear systems is a vast area with numerous approaches, e.g. \cite{meyn1993stability,pham2009contraction,bakry2013analysis,bolley2012convergence}. In order to derive convergence guarantees analogous to those available in linear systems, explicit quantitative  bounds are  required.

The motivation behind this paper is to derive quantitative convergence guarantees for stochastic adaptive control methods. To this end, we build upon methodologies at the intersection of stochastic differential equations and optimal transport \cite{eberle2019couplings,eberle2016reflection}. However, the existing methods in this area are too restrictive to be applied directly to common adaptive control schemes. In particular, these focus unconstrained processes with strong Lipschitz-like conditions on the drift term. However, in adaptive control, the parameters are typically constrained and their update rules often contain quadratic terms that violate the drift conditions. 

Our primary contribution to stochastic convergence theory is an extension of the methodology from \cite{eberle2019couplings} to constrained processes with less restrictive drift conditions. We derive an explicit exponential contraction bound under a specially constructed Wasserstein distance. 
The result implies exponential convergence  to a unique stationary distribution under a variety of measures, including total variation distance and Euclidean Wasserstein distances.
We then show how this result can be used to prove exponential convergence in a feedback-linearizable stochastic adaptive regulation problem. Additionally, we show how a projection method based on reflected stochastic differential equations can be used to constrain the parameters to an arbitrary closed convex set. This provides a flexible alternative to handling constraints, which contrasts with more specialized projection operators commonly employed in adaptive control \cite{lavretsky2013robust,hovakimyan2010L1,khalil2017high}.

The remaining parts of the paper are organized as follows. Section~\ref{sec:notation} presents
preliminary notation. Section~\ref{sec:general} presents the main contraction results, while
Section~\ref{sec:MRAC} presents the application to stochastic adaptive control. Section~\ref{sec:numerical} presents numerical
results and we provide closing remarks in Section~\ref{sec:conclusion}.
The main contraction theorem is proved in the appendix.

\section{Notation}
\label{sec:notation}

Random variables are denoted in bold, e.g. $\xb$. Time indices are denoted by subscripts, e.g. $\xb_t$ denotes a stochastic process.
We equip $\R^n$ with an inner product and norm denoted by $\langle \cdot ,\cdot \rangle$ and 
$\| \cdot \|$ respectively. We interpret $x,y\in \R^n$  as column vectors and let $x^\star$
denote the dual row vector such that $\langle x,y \rangle = x^\star y$. (Since $\langle\cdot,
\cdot\rangle$ is not necessarily the Euclidean inner product, we may have $x^\star \ne x^\top$.)
More generally, for a matrix $G$, let $G^\star$ denote its conjugate with respect to the inner
product. If $A$ is a square matrix, its trace is denoted by $\tr(A)$. 

If $\Xc$ is a closed convex set in $\R^n$, then at any $y\in\R^n$ the normal cone of $\Xc$ is defined as
\begin{equation}\label{normalCone}
N_{\Xc}(y) = \{v\in\R^n \ | \langle x-y, v \rangle \leq 0 \ \ \forall x\in\Xc \},
\end{equation}
and the convex projection on $\Xc$ is $\Pi_{\Xc}(y) = \arg\min_{x\in\Xc}\limits\|y-x\|$.

We use the shorthand notations $a\land b = \min\{a,b\}$ and $a\lor b = \max\{a,b\}$. $\Ind$ denotes the indicator function, and $I_d$ is the $d\times d$ identity matrix.

\section{Contraction for Reflected Stochastic Differential Equations}
\label{sec:general}

\begin{figure}
  \centering
  \includegraphics[width=.6\columnwidth]{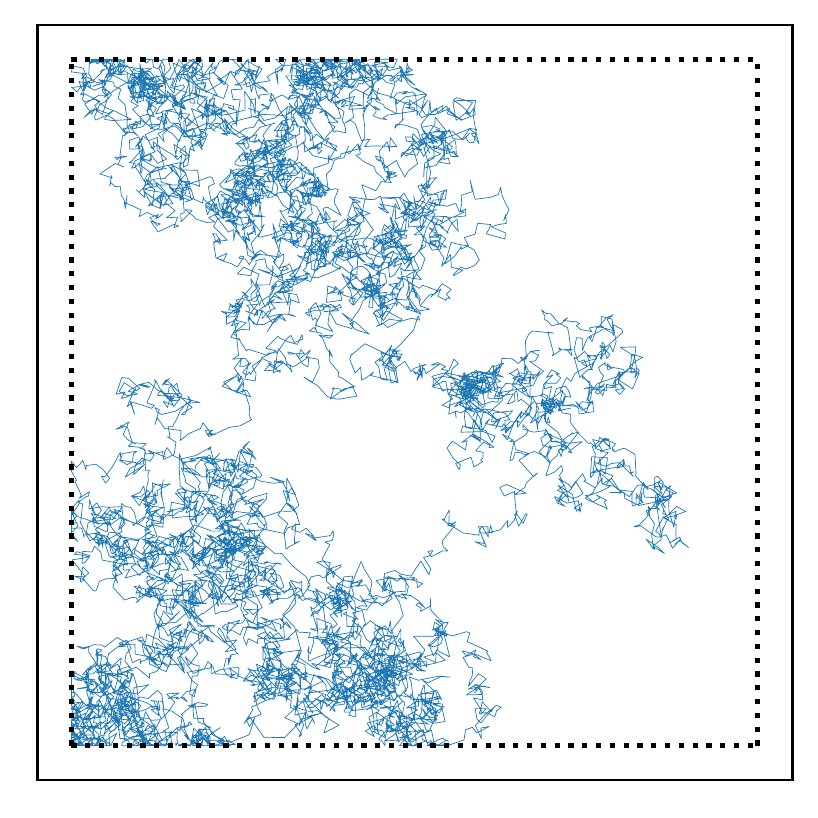}
  \caption{\label{fig:reflectedBM} {\bf Reflected Brownian Motion.} A non-smooth reflection process forces the trajectory to remain in a constraint set. }
\end{figure}

This section gives a general convergence result for reflected stochastic differential equations (RSDEs) over closed convex domains. See Fig.~\ref{fig:reflectedBM}.
In this paper, we consider RSDEs to handle the constraints that arise in adaptive parameter tuning rules.
A reflected stochastic differential equation is a stochastic differential equation which has been augmented with a special process which ensures that the trajectory remains within a constraint set.  The results of this section are required to handle the typical constraints placed on parameters in adaptive control methods. 

The results build upon the unconstrained contraction theory of \cite{eberle2019quantitative}, but substantial novel work is required to enable the  adaptive control applications in Section~\ref{sec:MRAC}.
In particular, we examine reflected SDEs to handle constraints, while \cite{eberle2019quantitative} considers unconstrained SDEs.

\subsection{Problem Setup}

Let $\Xc$ be a closed convex subset of $\R^n$. We will examine contractivity properties of reflected stochastic differential equations of the form:
\begin{equation}\label{diffX}
\d\xb_t = \Hs(\xb_t)\dt + G\d\wb_t - \vb_t\d\mub(t) ,
\end{equation}
where $G$ is an invertible $n\times n$ matrix, $\wb_t$ is a standard Brownian motion, and $\psib_t=-\int_0^t \vb_s\d\mub(s)$ is a bounded variation reflection process that enforces that $\xb_t \in \Xc$ for all $t\ge 0$ whenever $\xb_0\in\Xc$. In this case, $\wb_t$ has mean zero and $\E[\wb_t \wb_t^\star]=tI_n$. 

When $\Hs$ is Lipschitz, it can be shown that $\psib_t$ is the unique bounded variation process such that for all $t\ge 0$ and $\vb_t \in N_{\Xc}(\xb_t)$, $\|\vb_t\|\in \{0,1\}$,  $\mub$ is a random measure with $\mub([0,t])<\infty$,  and the solution has $\xb_t\in \Xc$. See \cite{tanaka1979stochastic,lions1984stochastic}. The Lipschitz condition can be relaxed to $\Hs$ being locally Lipschitz, provided that the process is non-explosive. See Section 2.4 of \cite{pilipenko2014introduction}.

Reflected stochastic differential equations can be simulated numerically via a projected Euler method:
$$
\xb_{t+\eta} \approx \Pi_{\Xc}(\xb_t+\eta \Hs(\xb_t)  + G (\wb_{t+\eta}-\wb_t)).
$$
See \cite{tanaka1979stochastic,slominski2001euler}. In other words, the effect of $\psib_t$ is the continuous time limit of a projection operation. 

\subsection{Assumptions}

The first requirement is a  one-sided growth condition. We assume that there is a function $\kappa(r):(0,\infty)\to [0,\infty)$ with $\int_0^1 s\kappa(s)ds <\infty$ and a non-negative number $\alpha$ such that for all $x,y\in \Xc$ with $r=\|x-y\|$ the following bound holds: 
\begin{equation}
  \label{eq:oneSided}
  \langle x-y, \Hs(x)-\Hs(y) \rangle \le \kappa(r)r^2 + \alpha r(\|x\|+\|y\|)
\end{equation}
This one-sided growth condition generalizes the one-sided Lipschitz condition from \cite{eberle2019quantitative}, which corresponds to the special case with $\alpha=0$. The extra terms are required to handle the application to adaptive control in Section~\ref{sec:MRAC}.

Let $\Ac$ denote the generator associated with the process $\xb_t$. Specifically, for any function $g:\Xc\to \R$
$$
(\Ac g)(x) = \lim_{h\downarrow 0} h^{-1}(\E[g(\xb_{h})|\xb_0=x]-g(x)).
$$

We assume that there is a twice continuously differentiable Lyapunov function $\Vs:\Xc\to [0,\infty)$ and positive numbers $\lambda$ and $C$ such that for all $x\in \Xc$ and all $t\ge 0$:
\begin{equation}
  \label{eq:fosterLyap}
(\Ac \Vs)(x) \le C -\lambda \Vs(x).
\end{equation}

We assume that $\Vs(x)$ increases monotonically with $\|x\|$. Specifically, there  is a strictly monotonically increasing function $\phi$ such that $\Vs(x)=\phi(\|x\|)$.  We will further assume that $\Vs$ grows at least linearly with $\|x\|$.

Let $R_1$ be the diameter of the set $\{(x,y)\in\Xc | \Vs(x)+\Vs(y) \le 4C/\lambda \}$. Linear growth implies that $R_1$ is finite. By construction, if $\|x-y\| > R_1$, then
\begin{equation}
\label{eq:linDecrease}
(\Ac\Vs)(x)+(\Ac\Vs)(y) \le -(\lambda/2)(\Vs(x)+\Vs(y)).
\end{equation}

Let $M$ be a positive number such that $M \ge R_1$ and for all $x$ with $\|x\|\ge M$, the following bound holds:
\begin{equation}
\label{eq:linGrowth}
\Vs(x) \ge \max\left\{ \frac{2}{\lambda} (\alpha\|x\| + C), \frac{4C}{\lambda}\left(2\|x\|+1 \right)\right\}.
\end{equation}
Let $R_2$ be the diameter of $\{(x,y)\in \Xc | \|x\| \le M \textrm{ and } \|y\|\le M \}$. Note that $R_2 \le 2M$ by the triangle inequality. 

In Section~\ref{sec:MRAC} we take $\Vs(x)=\|x\|^2+1$, so the growth conditions are automatically satisfied.

\subsection{Background on Wasserstein Distances}
Our main theory describes convergence in a Wasserstein distance.
To state the result, some basic concepts from optimal transport are required. See \cite{villani2008optimal} for a general reference. Let $P$ and $Q$ be probability measures over $\Xc$ with respect to the standard Borel sigma algebra. A measure, $\Gamma$, over $\Xc\times \Xc$ is called a \emph{coupling} of $P$ and $Q$ if its marginals are $P$ and $Q$, respectively. In other words, for any Borel set $S$, we have $\Gamma(S\times \Xc) = P(S)$ and $\Gamma(\Xc\times S) = Q(S)$. Let $\Cs(P,Q)$ denote the set of all couplings of $P$ and $Q$. 

If $\rho : \R^n \times \R^n \to [0,\infty)$ is a metric, the induced $q$-Wasserstein distance between $P$ and $Q$ is defined by:
\begin{equation}
\label{eq:wassersteinDef}
W^q_{\rho}(P,Q)=\inf_{\Gamma \in \Cs(P,Q)} \left(\int_{\Xc\times \Xc} \rho(x,y)^q\d\Gamma(x,y)\right)^{1/q}
\end{equation}

For simple notation, we follow the convention that $W_{\rho}:= W_\rho^1$ for general $\rho$ and for the norm,   $W^q := W^q_{\| \cdot \|}$. 

\subsection{Main Contraction Result}

The idea behind \cite{eberle2019quantitative} is to construct a new metric for which convergence of $W_{\rho}$ can be tractably analyzed, and then use the result to examine more standard measures such as $W^q$ and the total variation distance. We follow a similar approach, but the metric must be modified to account for the more general growth condition. 

Our metric will have the form:
\begin{multline}
  \rho(x,y) = \left[
    f(\|x-y\|) + \gamma \Vs(x)+\gamma \Vs(y) + \right. \\
  \Vs(x)\lor \phi(M) + \Vs(y) \lor \phi(M)
      \left.
    \right] \Ind(x\ne y).
  \end{multline}
  Recall that $\phi$ is a function such that $V(x)=\phi(\|x\|)$ and the $\Ind$ is the indicator function. Here $\gamma$ is a positive number defined below. 

  The function $f:[0,\infty) \to [0,\infty)$ is defined via the following chain of definitions:
  \begin{subequations}
    \label{eq:fDef}
    \begin{align}
      h(r) &= \frac{1}{\sigma_{\min}(G)^2}\left(\frac{1}{2} \int_0^r s\kappa(s)\d s + \alpha Mr \right)\\
      \varphi(r) &= e^{-h(r)} \\
      \Phi(r) &= \int_0^r \varphi(s)\d s \\
      \xi^{-1} &= \int_0^{R_1} \varphi(s)^{-1}\d s \\
      \beta^{-1} &= \int_0^{R_2} (\Phi(s)/\varphi(s))\d s \\
      g(r) &= 1-\frac{\xi}{4} \int_0^{r\land R_1} \varphi(s)^{-1}\d s -
             \frac{\beta}{4}\int_0^{r\land R_2}(\Phi(s)/\varphi(s))\d s \\
      \label{eq:fInt}
      f(r) &= \int_0^{r\land R_2} \varphi(s) g(s)\d s.
  \end{align}
  \end{subequations}
  Here $\sigma_{\min}(G)>0$ is the smallest singular value of $G$.
  
  Additionally, we set
  \begin{equation}
    \label{eq:gamma}
    \gamma = \frac{\xi \sigma_{\min}(G)^2}{4C}.
  \end{equation}

  The general contraction result is given below. The proof is given in appendix~\ref{appendix}.
  
\begin{theorem}
  \label{thm:genContraction}
  {\it 
    The function $\rho(x,y)$ 
 is a metric over $\Xc$. Let $\xb_t$ and $\yb_t$ be two solutions to (\ref{diffX}) with respective laws $\P_t$ and $\Q_t$.  If the initial laws satisfy $\int_{\Xc} \Vs(x)d\P_0(x)<\infty$ and $\int_{\Xc}\Vs(x)d\Q_0(x)$, then
 $$
 W_{\rho}(\P_t,\Q_t)\le e^{-at} W_{\rho}(\P_0,\Q_0).
 $$
 where $a = \min\{ \lambda, \xi \sigma_{\min}(G)^2 , \beta \sigma_{\min}(G)^2\}/2$.
}

The following corollary establishes convergence to a stationary distribution in total variation distance and norm-based $q$-Wassertein distances. It is analogous to Corollary 2.1 and Remark 3.4 of \cite{eberle2019quantitative}. The proof is omitted, since the argument from \cite{eberle2019quantitative} works with minimal modification. 

\begin{corollary}
  \label{cor:equilibrium}
  {\it
  The system from (\ref{diffX})
  has a unique stationary distribution $\pi$ such that $\int_{\Xc} \Vs(x) d\pi(x) < \infty$.
  Let $\xb_t$ be a solution to (\ref{diffX}) with law $\P_t$ such that $\int_{\Xc} \Vs(x) d\P_0(x) <\infty$. If $\Vs(x)\ge 1$ for all $x\in \Xc$, then $\P_t$ converges in total variation as:
  \begin{equation*}
    \|\P_t-\pi\|_{\TV} \le \gamma^{-1} e^{-at} W_{\rho}(\P_0,\pi).
  \end{equation*}

  Say $q> 1$ and $\frac{1}{p}+\frac{1}{q}=1$ . If $\Vs(x) \ge \|x\|^q$ for all $x\in\Xc$, then  $\P_t$ converges with respect to $W^q$ as:
  \begin{equation*}
    W^q(\P_t,\pi) \le 2^{1/p} \left( \gamma^{-1} W_{\rho}(\P_0,\pi) \right)^{1/q} e^{-at/q}.
  \end{equation*}
  }
\end{corollary}
\end{theorem}

\subsection{Discussion}

We describe the distinctions between our results and those of \cite{eberle2019quantitative}. The most obvious is that ours applies to processes reflected to remain in the set $\Xc$, while \cite{eberle2019quantitative} examines unconstrained SDEs. Additionally, our one-sided growth condition is used instead of a one-sided Lipschitz condition: 
\begin{equation*}
\langle x-y,\Hs(x)-\Hs(y)\rangle \le \kappa(r) r^2.
\end{equation*}
with the same definition of $\kappa$.
This condition, however, fails in very basic versions of the control problem from Section~\ref{sec:MRAC}, thus necessitating the our condition \eqref{eq:oneSided}.
Our main application utilizes quadratic Lyapunov functions, and so we restrict to the case of Lyapunov functions $\Vs(x)$ with quadratic growth. This is less general than \cite{eberle2019quantitative}, but leads to substantially simpler analysis. Additionally, this leads to a much simpler and more explicit metric.

\section{Application to Adaptive Regulation}
\label{sec:MRAC}

\subsection{Problem Setup}

We analyze a stochastic variation of a model reference adaptive control problem examined in Chapter 9 of \cite{lavretsky2013robust}.

The basic plant model has the form
\begin{equation}
  \label{eq:plant}
\d\xb_t = [\bar{A} \xb_t + B(\bar{\Omega}^\top \Psi(\xb_t) + \ub_t)]\dt + G_x\d\wb_t^x.
\end{equation}
Here $\bar{A}$ is an unknown state matrix, $B$ is a known input matrix, $\Psi$ is a known vector of feature functions (which could be linear or nonlinear), $\bar{\Omega}$ is an unknown matrix of parameters,  and $G_x$ is an unknown matrix scaling the Brownian motion. The state is $\xb_t$ and the inputs are $\ub_t$. The setup in \cite{lavretsky2013robust} also includes an unknown scaling factor on $B$, which we have omitted for simplicity.   

We focus on the problem of adaptive regulation, while \cite{lavretsky2013robust} examines tracking problems. In our setup, the matching assumption is that there is a known Hurwitz matrix, $A$, and an unknown feedback gain, $\bar{K}$, such that
$$
\bar{A}+B\bar{K}=A.
$$
Here $A$ is the state matrix for the \emph{reference system}.

It appears to be possible to extend the contraction theory to  tracking problems, but this is left for future work.

If we knew $\bar{K}$ and $\bar{\Omega}$, we could set $\ub_t = \bar{K}\xb_t - \bar\Omega^\top \Psi(\xb_t)$ and render the system stochastically stable:
$$
\d\xb_t = A\xb_t\dt + G_x\d\wb_t^x.
$$
The challenge is that we do not know $\bar{K}$ or $\bar\Omega$. So instead, we use $\ub_t = \Kb_t \xb_t - \Omegab_t^\top \Psi(\xb_t)$, where $\Kb_t$ and $\Omegab_t$ are estimates. We derive rules for their computation later.

To simplify notation, we set
\begin{equation*}
  \bar\Theta = \begin{bmatrix}-\bar K^\top \\ \bar\Omega \end{bmatrix}, \:
  \Thetab_t = \begin{bmatrix}-\Kb_t^\top \\ \Omegab_t \end{bmatrix}, \:
  \Lambda(\xb_t) = \begin{bmatrix}
    \xb_t \\
    \Psi(\xb_t)
  \end{bmatrix}.
\end{equation*}
Then the dynamics of \eqref{eq:plant} with $\ub_t = \Kb_t \xb_t - \Omegab_t^\top \Psi(\xb_t)$ can be written as 
\begin{align}\label{xDot}
\d\xb_t = \big(A\,\xb_t +  B\,(\bar\Theta-\Thetab_t)^\top\Lambda(\xb_t)\big)\dt + \sigx\,\d\wb_t^x.
\end{align} 
 where fixed matrices $A\in\R^{n\times n}$ and $B\in\R^{n\times\ell}$ form a controllable pair
 and $\wb_t^x\in\R^n$ is standard Brownian motion with coefficient matrix $\sigx\in\R^{n\times
 n}$, which is fixed and invertible.

We assume that the function $\Lambda:\R^n\to\R^L$ is Lipschitz:
$$
\|\Lambda(x)-\Lambda(y)\|_2 \le \Lc \|x-y\|_2,
$$
for some $\Lc > 0$, where $\|\cdot \|_2$ is the Euclidean norm.

Assume that $\bar\Theta$ and $\Thetab_t$ are $L\times \ell$ matrices, and set $m=L\ell$. Let $S:\R^m \to \R^{L\times \ell}$ be the reshaping function defined by $S(v)_{i,j}=v_{(i-1)\ell+j}$ 
for $i=1,\ldots,L$ and $j=1,\ldots,\ell$. Then $S$ is an invertible linear function.

Let $\bar\theta=S^{-1}(\bar\Theta)$ be the unknown parameters and let $\Kc$ be a compact convex subset of $\R^m$, containing 
$\bar\theta$. We assume that $\Kc$ is known. Let $D$ be the diameter of $\Kc$.

Now let $\Xc=\R^n\times\Kc$ be the closed convex subset of $\Rnm$ containing the combined
state $\zb_t=[\xb_t^\top\ \thetab_t^\top]^\top$, and assume that $\thetab_t$ has dynamics of the form
\begin{align}\label{thhDiff}
\d\thetab_t = \Rs(\zb_t)\,\dt + G_\theta\,\d\wb_t^\theta -
\vb_t^\theta\,\d\mub^\theta(t) \ ,
\end{align}
where $\Rs:\Xc\to\R^m$, $\wb_t^\theta\in\R^m$ is standard Brownian motion
with invertible coefficient matrix $G_\theta\in\R^{m\times m}$, and
$\psib_t^\theta = -\int_0^t\vb_s^\theta\,\d\mub^\theta(s)$ is a bounded variation reflection process that
enforces $\thetab_t\in\Kc$ for all $t\geq0$ whenever $\thetab_0\in\Kc$. Here we assume that the reflections are computed with respect to the standard Euclidean inner product for simplicity.

Under any norm of the form $\|z\|^2 = x^\top Px + \theta^\top \theta$, the joint dynamics defined by \eqref{xDot} and \eqref{thhDiff} define a special case of the reflected stochastic differential equation \eqref{diffX}. Indeed, if $z^\top = [x^\top,\theta^\top]$, then every element of $N_{\Xc}(z)$ has the form $v^\top = [0^\top ,(v^\theta)^\top]$, where $v^\theta \in N_{\Kc}(\theta)$ and $N_{\Kc}(\theta)$ is defined with respect to the Euclidean inner product.

\begin{remark}
Our parameter tuning rule from \eqref{thhDiff} differs from typical methods from adaptive control in how it forces $\thetab_t$ to remain in the constraint set, $\Kc$. Indeed, our method uses reflection processes, which can be approximated in discrete time by convex projections.
For simple sets such as Euclidean balls and boxes, the convex projections have simple analytic formulas. More generally, for any convex set, the convex projection can be computed via optimization. In contrast, the parameters of adaptive control laws are commonly constrained using specialized projection operators that are designed for specific classes of convex sets \cite{lavretsky2013robust,hovakimyan2010L1,ioannou2006adaptive,khalil2017high}.
\end{remark}

\subsection{Lyapunov-Based Adaptation}
Here we follow a Lyapunov-based design procedure to design the update rule,  \eqref{thhDiff}, similar to the method from \cite{lavretsky2013robust}.
The main differences are that we examine stochastic problems and the constraints are enforced by reflection. 

First we construct the Lyapunov function candidate. Fix a positive definite $Q\in \R^{n\times n}$. Then since $A$ is Hurwitz, there is a unique positive definite $P$ such that
$$
A^\top P+P\,A=-Q.
$$

For $z^\top = [x^\top,\theta^\top]$,
we define the Lyapunov function candidate by:
\begin{equation}
\label{lyap}
\Vs(z) = x^\top Px +
(\theta-\bar\theta)^\top(\theta-\bar\theta) + 1.
\end{equation}
Define a norm over $\Rnm$ by $\|z\|^2 = x^\top Px + \theta^\top \theta$.

Theorem~\ref{thm:genContraction} 
requires that the Lyapunov function be a monotonic function of a norm. This could be attained using the affine coordinate transformation $\hat z^\top = [x^\top,(\theta-\bar\theta)^\top]$. In that case we can set $\hat\Vs(\hat z)=\Vs(z) = \|\hat z\|^2 +1$. In the analysis below, we will work in the original coordinates.

To derive the required decrease condition, we use It\^o's formula:
\begin{align}\nonumber
\d\Vs(\zb_t)  &=  \nabla_x\Vs^\top\d\xb_t 
 + \frac{1}{2}\,\d\xb_t^\top(\nabla_x^2\,\Vs)\,\d\xb_t \\ \nonumber
 &\hspace{11pt}+ \nabla_\theta\Vs^\top\d\thetab_t
 + \frac{1}{2}\,\d\thetab_t^\top(\nabla_\theta^2\,\Vs)\,\d\thetab_t \\
 \nonumber
&= 2\,\xb_t^\top P^\top\d\xb_t + \d\xb_t^\top P\,\d\xb_t \\ \nonumber
&\hspace{11pt}+ 2\,(\thetab_t-\bar\theta)^\top\d\thetab_t +
\d\thetab_t^\top\d\thetab_t \\
\nonumber
&= \xb_t^\top\big(A^\top P+PA\big)\,\xb_t\,\dt \\ \nonumber
&\hspace{11pt}- 2\,\xb_t^\top P\,B(\Thetab_t-\bar\Theta)^\top\Lambda(\xb_t)\,\dt 
+ 2\,\xb_t^\top P\,\sigx\,\d\wb_t^x \\ \nonumber
&\hspace{11pt}+ \d\xb_t^\top P\,\d\xb_t 
+ 2\,\tr\big((\Thetab_t-\bar\Theta)^\top \d\Thetab_t\big)
+\d\thetab_t^\top\d\thetab_t  \\ \nonumber
&= -\xb_t^\top Q\,\xb_t\,\dt \\ \nonumber
&\hspace{11pt}+ 2\,\tr\big((\Thetab_t-\bar\Theta)^\top(\d\Thetab_t
-\Lambda(\xb_t)\,\xb_t^\top P\,B\,\dt)\big) \\
\label{lyapDiff}
&\hspace{11pt} + \tr(\sigx^\top P\,\sigx)\,\dt + \d\thetab_t^\top\d\thetahb_t
+ 2\,\xb_t^\top P\,\sigx\,\d\wb_t^x \ ,
\end{align}
where the third equality uses symmetry of $P$ and that
 \begin{equation*}
 2\,\xb_t^\top P^\top\d\xb_t = \xb_t^\top P^\top\d\xb_t + \d\xb_t^\top P\,\xb_t \ ,
 \end{equation*}
 and 
the fourth equality uses the definition of $Q,P$. Additionally, is uses the fact that
\begin{equation*}
\xb_t^\top P\,B(\Thetab_t-\bar\Theta)^\top\Lambda(\xb_t) =
\tr\big((\Thetab_t-\bar\Theta)^\top\Lambda(\xb_t)\,\xb_t^\top P\,B\big) ,
\end{equation*}
and the It\^o rule $\d\wb_t^x(\d\wb_t^x)^\top = \dt\,I_n$.

By examining the second term in \eqref{lyapDiff}, we see that the $\Lambda(\xb_t)\,\xb_t^\top
P\,B\,\dt$ term is canceled if $\Rs$ in \eqref{thhDiff} is defined as
\begin{align}\label{Rdef}
\Rs(\zb_t) = S^{-1}\big(\Lambda(\xb_t)\,\xb_t^\top P\,B\big) \ ,
\end{align}
where $S^{-1}$ is the inverse shaping function.

Plugging \eqref{thhDiff} into \eqref{lyapDiff}, as $\d\Thetab_t=S(\d\thetab_t)$, then gives
\begin{align}\nonumber
&\d\Vs = -\xb_t^\top Q\xb_t\,\dt \\ \nonumber
&\hspace{11pt}+ 2\,\tr\big((\Thetab_t-\bar\Theta)^\top \big(S(G_\theta\,\d\wb_t^\theta) -
S(\vb_t^\theta\,\d\mub^\theta(t))\big)\big) \\ \nonumber
&\hspace{11pt} + \tr(\sigx^\top P\,\sigx)\,\dt + \tr(G_\theta^\top G_\theta)\,\dt
+ 2\,\xb_t^\top P\,\sigx\,\d\wb_t^x \\
\nonumber
&= \big(-\xb_t^\top Q\xb_t + \tr(\sigx^\top P\,\sigx) + \tr(G_\theta^\top G_\theta)\,\big)\,\dt
\\ \nonumber
&\hspace{11pt} + 2\,\xb_t^\top P\,\sigx\,\d\wb_t^x +
2\,(\thetab_t-\bar\theta)^\top G_\theta\,\d\wb_t^\theta \\ \label{lyapDiff1}
&\hspace{11pt}- 2\,(\thetab_t-\bar\theta)^\top\vb_t^\theta\,\d\mub^\theta(t) \ ,
\end{align}
where again we use that $\d\wb_t^\theta (\d\wb_t^\theta)^\top = \dt\,I_m$. 

Now since
$\vb_t^\theta\in N_{\Kc}(\thetab_t)$  and
$\mub^\thetah$ is a nonnegative measure,  \eqref{normalCone} implies
\begin{align*}
- (\thetab_t-\bar\theta)^\top\vb_t^\theta\,\d\mub^\theta(t) \leq 0 \ .
\end{align*}

Therefore, the  generator of the Lyapunov function $\Vs(\zb_t)$ satisfies
\begin{equation}
\label{eq:adaptiveGenerator}
\Ac\Vs(z) \le \big(-x^\top Qx + \tr(\sigx^\top P\,\sigx) + \tr(G_\theta^\top G_\theta)\,\big).
\end{equation}
Now using standard quadratic form bounds, followed by the diameter condition on $\Kc$ gives:
\begin{align*}
  x^\top Q x &\ge \frac{\lambda_{\min}(Q)}{\lambda_{\max}(P)} x^\top P x \\
  &\ge \frac{\lambda_{\min}(Q)}{\lambda_{\max}(P)} \Vs(z) - \frac{\lambda_{\min}(Q)}{\lambda_{\max}(P)}(D^2+1) \ .
\end{align*}
Here $\lambda_{\min}(Q)$ is the minimum eigenvalue of $Q$ and $\lambda_{\max}(P)$ is the maximum eigenvalue of $P$. Plugging this bound into \eqref{eq:adaptiveGenerator} shows $\Vs$ satisfies \eqref{eq:fosterLyap} with
\begin{align*}
  C &= \tr(\sigx^\top P\,\sigx) + \tr(G_\theta^\top G_\theta) +   \frac{\lambda_{\min}(Q)}{\lambda_{\max}(P)}(D^2+1)  \\
  \lambda &=  \frac{\lambda_{\min}(Q)}{\lambda_{\max}(P)}.
\end{align*}

In particular, $\Vs$ satisfies all of the required conditions to apply the general theory from Section~\ref{sec:general}.

\subsection{One-Sided Growth for Adaptive Regulation}

The final task needed to apply Theorem~\ref{thm:genContraction} to the adaptive regulation problem is ensuring that the one-sided growth condition from \eqref{eq:oneSided} holds. Note that the combination of \eqref{xDot}, \eqref{thhDiff}, \eqref{Rdef} leads to a special case of \eqref{diffX} with:

\begin{align}\nonumber
&\Hs(\zb_t) =
\begin{bmatrix}
A\,\xb_t  +B\,\big(\bar\Theta-\Thetab_t)^\top\Lambda(\xb_t)\big) \\[6pt]
S^{-1}\big(\Lambda(\xb_t)\xb_t^\top PB\big)
\end{bmatrix} \\ \nonumber
&G =
\begin{bmatrix}
\sigx & 0 \\
0 & G_\theta
\end{bmatrix} , \ \ 
\vb_t\,\d\mub_t =
\begin{bmatrix}
0 \\
\vb_t^\theta\,\d\mub(t)^\theta
\end{bmatrix}.
\end{align}

Set $z^\top = [x^\top,\theta^\top]$ and $\tilde z^\top = [\tilde x^\top,\tilde \theta^\top]$. Direct calculation using the specialized inner product shows that 
\begin{align}
  \nonumber
  & 
\langle z-\tilde z,\Hs(z)-\Hs(\tilde z) \rangle \\
  \nonumber
&= (x-\xw)^\top P\,A\,(x-\xw) \\
  \nonumber
&\hspace{11pt} -(x-\xw)^\top P\,B\,(\Theta-\bar\Theta)^\top\Lambda(x) \\
  \nonumber
&\hspace{11pt}  + (x-\xw)^\top P\,B\,(\Thetaw-\bar\Theta)^\top\Lambda(\xw)\\
  \nonumber
&\hspace{11pt} +\tr\Big((\Theta-\Thetaw)^\top(\Lambda(x)x^\top-
\Lambda(\xw)\xw^\top\big)\,PB\Big) \\
\nonumber
&=
-\frac{1}{2}(x-\xw)^\top Q_x\,(x-\xw) \\
\nonumber
&\hspace{11pt}
-\tr\Big((\Theta-\bar\Theta)^\top\Lambda(x)(x-\xw)^\top P\,B\Big) \\
\nonumber
&\hspace{11pt} + \tr\Big((\Thetaw-\bar\Theta)^\top\Lambda(\xw)(x-\xw)^\top P\,B\Big)\\
  \nonumber
&\hspace{11pt} +\tr\Big((\Theta-\Thetaw)^\top(\Lambda(x)x^\top-
\Lambda(\xw)\xw^\top\big)\,PB\Big)
\end{align}
\begin{align}
\nonumber
&= -\frac{1}{2}(x-\xw)^\top Q\,(x-\xw) \\ \nonumber
&\hspace{11pt}-\cancel{\tr\Big(\Theta^\top\Lambda(x)x^\top P\,B\Big)}
+\tr\Big(\Theta^\top\Lambda(x)\xw^\top P\,B\Big) \\
\nonumber
&\hspace{11pt} +\tr\Big(\bar\Theta^\top\Lambda(x)x^\top P\,B\Big)
-\tr\Big(\bar\Theta^\top\Lambda(x)\xw^\top P\,B\Big) \\
\nonumber
&\hspace{11pt} +\tr\Big(\Thetaw^\top\Lambda(\xw)x^\top P\,B\Big)
-\cancel{\tr\Big(\Thetaw^\top\Lambda(\xw)\xw^\top P\,B\Big)} \\
\nonumber
&\hspace{11pt} -\tr\Big(\bar\Theta^\top\Lambda(\xw)x^\top P\,B\Big)
+\tr\Big(\bar\Theta^\top\Lambda(\xw)\xw^\top P\,B\Big) \\
\nonumber
&\hspace{11pt} +\cancel{\tr\Big(\Theta^\top\Lambda(x)x^\top P\,B\Big)}
-\tr\Big(\Theta^\top\Lambda(\xw)\xw^\top P\,B\Big) \\
\nonumber
&\hspace{11pt} -\tr\Big(\Thetaw^\top\Lambda(x)x^\top P\,B\Big)
+\cancel{\tr\Big(\Thetaw^\top\Lambda(\xw)\xw^\top P\,B\Big)} \\ \nonumber
&\leq \hspace{9.75pt} \tr\Big(\Theta^\top\Lambda(x)\xw^\top P\,B\Big)
+\tr\Big(\Thetabar^\top\Lambda(x)x^\top P\,B\Big) \\ \nonumber
&\hspace{11pt}-\tr\Big(\Thetabar^\top\Lambda(x)\xw^\top P\,B\Big)
+\tr\Big(\Thetaw^\top\Lambda(\xw)x^\top P\,B\Big) \\ \nonumber
&\hspace{11pt}+\tr\Big(\Thetabar^\top\Lambda(\xw)\xw^\top P\,B\Big)
-\tr\Big(\Thetabar^\top\Lambda(\xw)x^\top P\,B\Big) \\ \nonumber
&\hspace{11pt}-\tr\Big(\Theta^\top\Lambda(\xw)\xw^\top P\,B\Big)
-\tr\Big(\Thetaw^\top\Lambda(x)x^\top P\,B\Big) \\
  \label{eq:cancelations}
  &=
 - x^\top P\,B\,(\Thetaw-\bar\Theta)^\top(\Lambda(x)-\Lambda(\xw)) \\
  \nonumber
& \hspace{11pt} + \xw^\top P\,B\,(\Theta-\bar\Theta)^\top(\Lambda(x)-\Lambda(\xw)),
\end{align}
where we use the fact that $a^\top b = \tr(b\,a^\top)$ holds for any two equal length vectors $a,b$, and in the last inequality we drop the nonpositive $-\frac{1}{2}(x-\xw)^\top Q\,(x-\xw)$ term.\\
Set
\begin{align*}
  y &= P\,B\,(\Theta-\bar\Theta)^\top(\Lambda(x)-\Lambda(\xw)) \\
  \tilde y &= P\,B\,(\Thetaw-\bar\Theta)^\top(\Lambda(x)-\Lambda(\xw)).
\end{align*}
Then note that
\begin{align*}
  \|y-\tilde y\|_2 &= \| PB(\Theta - \Thetaw)^\top (\Lambda(x)-\Lambda(\tilde x))\|_2 \\
                   &\le \|PB\|_2 \|\Theta - \Thetaw\|_2 \,\Lc \|x-\tilde x\|_2 \\
                   &\le \|PB\|_2 \|\Theta - \Thetaw\|_F \,\Lc \|x-\tilde x\|_2 \\
  &\le \|PB\|_2 \,D\,\Lc \|x-\tilde x\|_2 \ .
\end{align*}
Here $\|\cdot \|_2$ applied to matrices refers to the induced $2$-norm. Then the first inequality uses submultiplicativity followed by the Lipschitz assumption on $\Lambda$. Next we note that the induced $2$-norm is bounded above by the Frobenius norm, and that $\|\Theta-\Thetaw\|_F = \|\theta-\tilde\theta\|_2$. So, the final inequality follows from the diameter bound.\\\\
An analogous derivation shows that 
\begin{align*}
  \|\tilde y\|_2 &\le \|PB\|_2 \,D\,\Lc \|x-\tilde x\|_2.
\end{align*}
Now the right side of \eqref{eq:cancelations} can be upper bounded by:
\begin{align*}
  -x^\top \tilde y + \tilde x^\top y &= (\tilde x-x)^\top \tilde y +\tilde x^\top (y-\tilde y) \\
                                     &\le \|\tilde x-x\|_2 \|\tilde y\|_2 + \|\tilde x\|_2 \|y-\tilde y\|_2 \\
                                     &\le \|PB\|_2 \,D\,\Lc\left(
                                       \|x-\tilde x\|_2^2 + \|\tilde x\|_2 \| x-\tilde x\|_2
                                       \right).
\end{align*}
Now using the bound $\|x\|_2 \le \frac{1}{\sqrt{\lambda_{\min}(P)}} \|z\|$ gives a special case of \eqref{eq:oneSided} with
$
  \kappa(r)=\alpha = \frac{ \|PB\|_2 \,D\,\Lc}{\lambda_{\min}(P)}.
$

\newpage
The results are summarized in the following theorem:

\begin{theorem}
  {\it
  The controller $\ub_t = \Thetab_t^\top \Psi(\xb_t)$ with parameters by \eqref{thhDiff} and \eqref{Rdef} drives the closed-loop system to a unique stationary distribution. Convergence is exponential with respect to total variation distance and $W^2$, with convergence rates described by Theorem~\ref{thm:genContraction} and Corollary~\ref{cor:equilibrium}.
  }
\end{theorem}

\section{Numerical Results}
\label{sec:numerical}

\begin{figure}[h]
\centering
\includegraphics[width=.9\columnwidth]{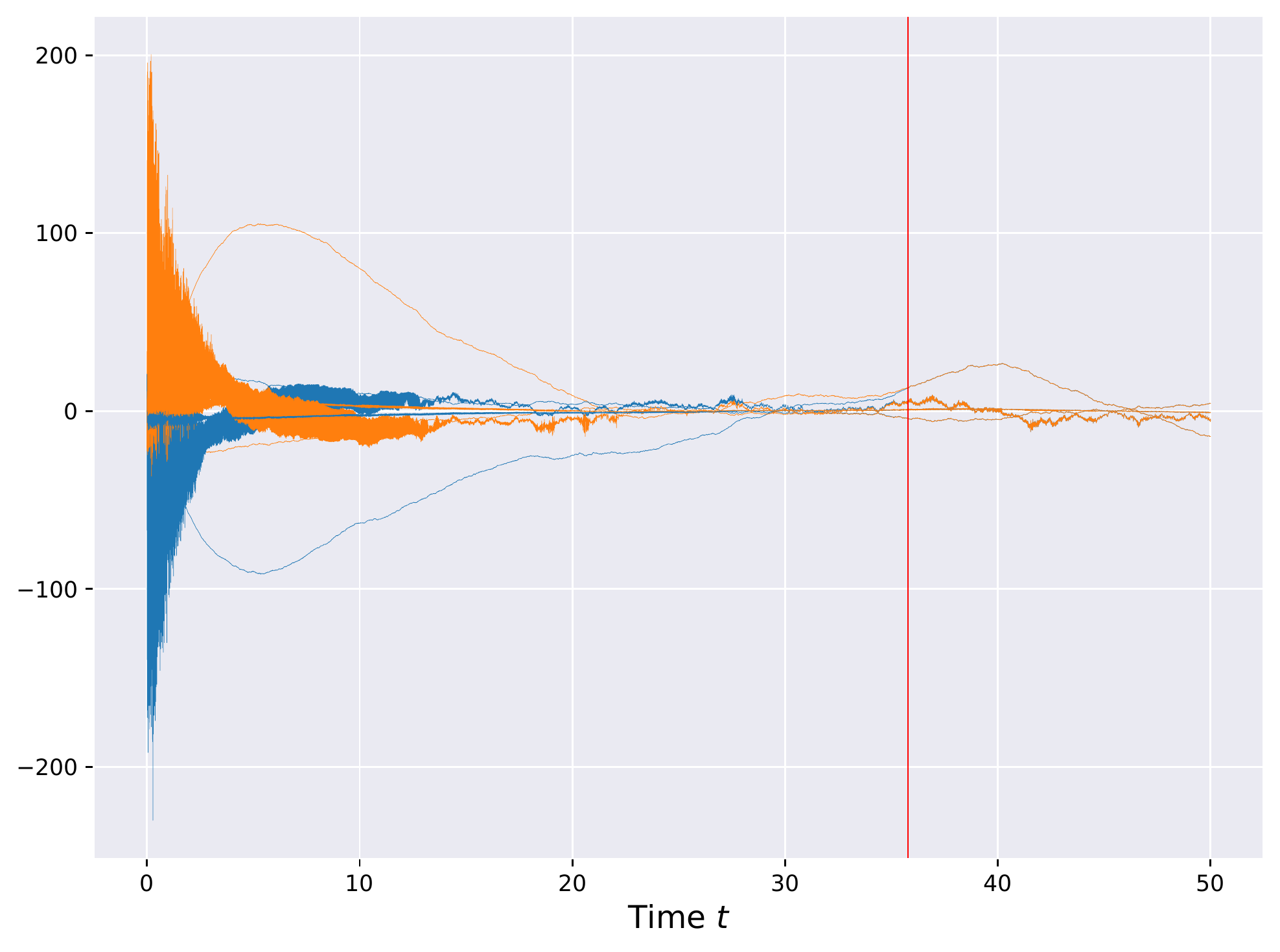}
\caption{\label{fig:coupledX} {\bf Reflection Coupling.} System states $\xb_t$ and $\xbw_t$ of the
	reflection coupling. The red line shows the coupling time $\taub$. }
\end{figure}

\begin{figure}[h]
	\centering
	\includegraphics[width=.9\columnwidth]{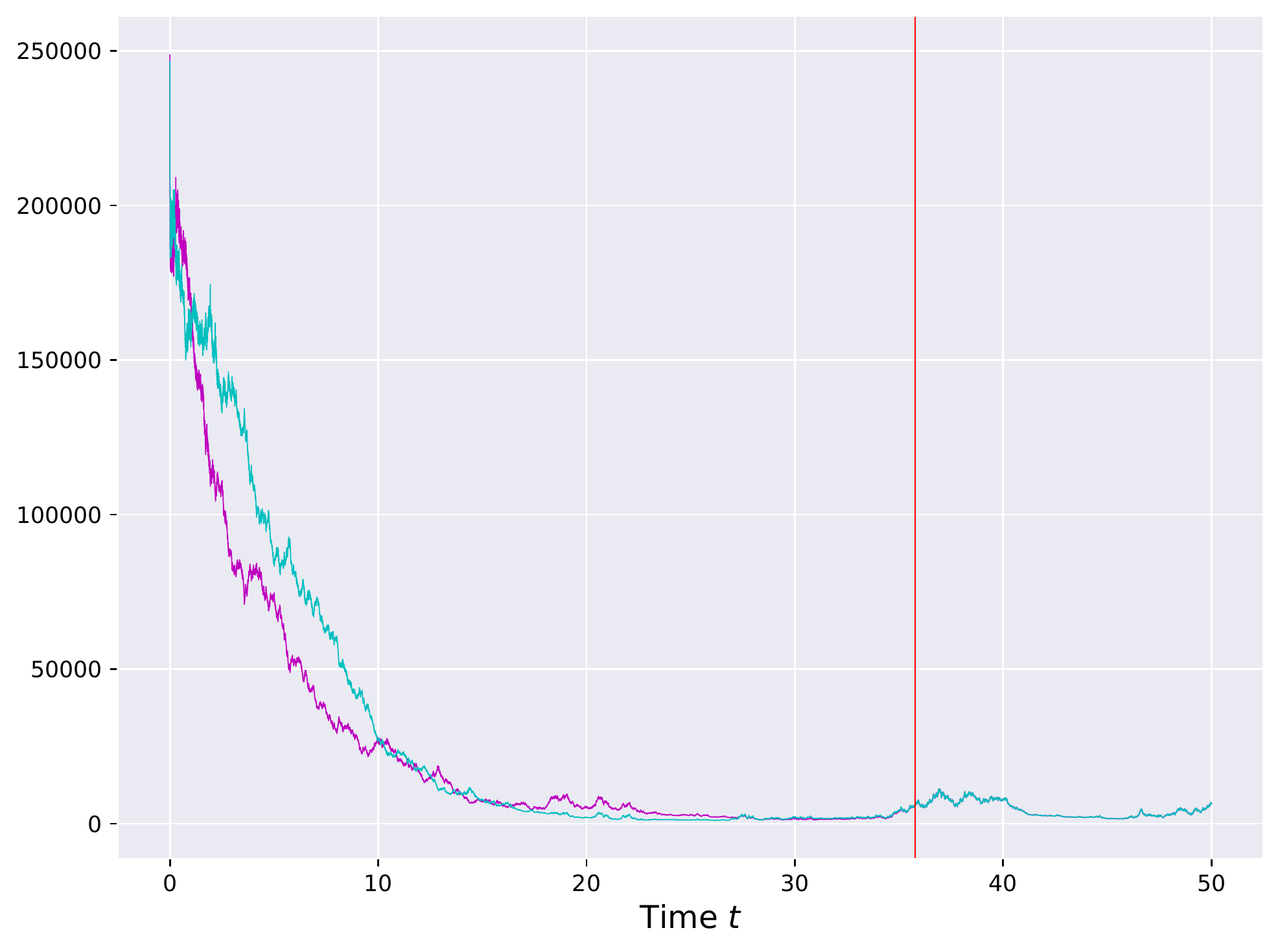}
	\caption{\label{fig:lyaps} {\bf Lyapunov Functions.} Lyapunov functions $\Vs(\zb_t)$ and $\Vs(\zbw_t)$ of the reflection coupling. The red line shows the coupling time $\taub$. }
\end{figure}

\begin{figure}[h]
	\centering
	\includegraphics[width=.9\columnwidth]{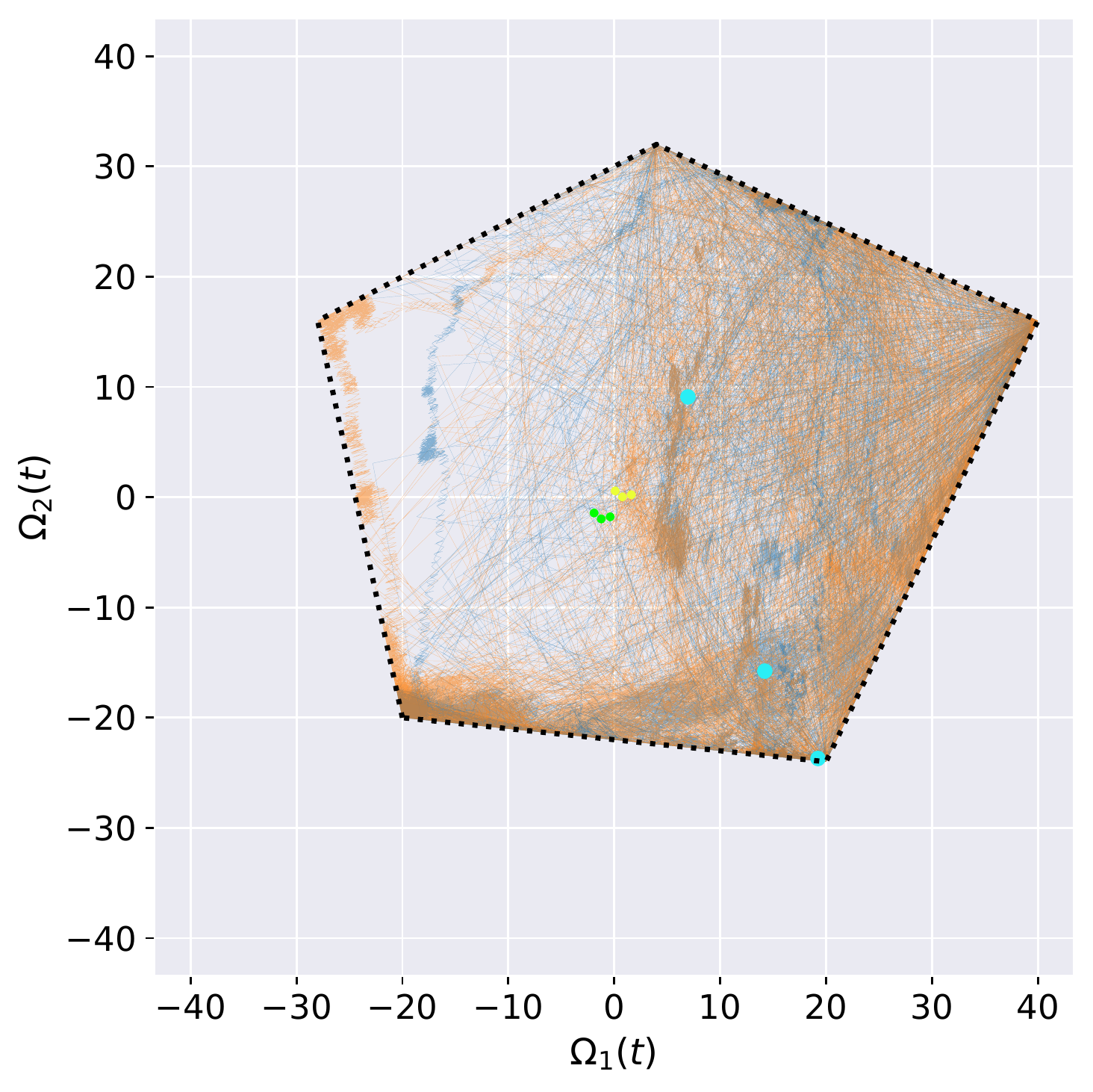}
	\caption{\label{fig:proj} {\bf Convex Projection.} Projection of each row of parameter
	estimates $\Omegab_t$ and $\widetilde\Omegab_t$ of the reflection coupling on the same
	polygon. Green and yellow dots show initial, while cyan dots show final (coupled) points. }
\end{figure}

We simulate\footnote{All code available at: https://github.com/tylerlekang/CDC2021} the plant system from Sec. 11.5 of \cite{lavretsky2013robust}, with $n=4$, $\ell=2$,
$N=3$, $L=7$, and thus $m=14$ (flattening of $\bar\Omega$ and $\bar K$). The $G$ matrices were $G_x = I_n$ and $G_\theta=I_m$. The Euler method was run over time $t\in[0,50]$ seconds
and with timestep $\eta=0.001$ seconds. The compact set $\Kc$ was a separate 2 dimensional
polygon applied to each row of $\Thetab_t$ (a parameter for each control input). For the $N=3$ rows of $\Omegab_t$ ($6$ total parameters), the same polygon was used, and is shown in Figure~\ref{fig:proj}. Figure~\ref{fig:coupledX} shows the reflection coupling of the system
states $\xb_t$ and $\xbw_t$, and their coupling time $\taub$. Figure~\ref{fig:lyaps} shows the
respective Lyapunov functions $\Vs(\zb_t)$ and $\Vs(\zbw_t)$, also with the coupling time.
Figure~\ref{fig:proj} shows the results of the convex projection on the polygon for each 2-dimensional space corresponding to a row of $\Omegab_t$ and $\widetilde\Omegab_t$.

\section{Conclusion}
\label{sec:conclusion}

\subsection{Discussion on Practical Applications}

In regards to practical application of the results, the authors would like to highlight
two key factors: 1) the flexibility afforded by the projection method in the various geometries
that can constrain the parameter estimates, as an improvement over existing methods, and 2)
opening the door for analysis of Langevin Algorithms on more general state spaces (see
\cite{lamperski2021projected}).

\subsection{Closing Remarks}

In this paper we introduced a novel extension of contraction methods for SDEs which enables restrictions to closed convex domains. We
utilized this theory to prove convergence for stochastic versions of adaptive controllers from \cite{lavretsky2013robust}. Future
work includes expanding the class of systems for which this theory holds, including exogenous
input tracking.

\printbibliography

\appendices
\section{Proof of Theorem~\ref{thm:genContraction}}\label{appendix}
Many of the arguments are similar to those of \cite{eberle2019quantitative}. We give terse explanations of these, and then highlight the differences.

The main idea is to define an explicit coupling between $\xb_t$ and $\yb_t$ and  show that $e^{at} \rho(\xb_t,\yb_t)$ is a supermartingale under this coupling. Then the definition of $W_\rho$ from \eqref{eq:wassersteinDef} followed by the supermartingale property show that:
$$
W_\rho(\P_t,\Q_t) \le \E[\rho(\xb_t,\yb_t)] \le e^{-at} \E[\rho(\xb_0,\yb_0)].
$$
The final bound from the theorem follows by taking an optimal coupling for $\P_0$ and $\Q_0$. This optimal coupling exists by Theorem 4.1 of \cite{villani2008optimal}.

\subsection*{Showing That $\rho$ is a Metric}
First note that
$$
[\gamma V(x) + \gamma V(y) + V(x)\lor \phi(M) + V(y)\lor\phi(M)] \Ind(x\ne y)
$$
is a metric since $V(x)\ge 0$ and $\phi(M)>0$. In particular, it is a weighted version of the Hamming metric.

Now since the sum of metrics is again a metric, it suffices to show that $f(\|x-y\|)$ is a metric. This holds provided that $f$ is concave, with $f'(0)>0$ and $f(0) =0$.  See
 \cite{eberle2019quantitative,eberle2016reflection} for details.

To this end, it suffices to show that $f(0)=0$, $f$ is monotonically increasing, and $f(a+b)\le f(a)+f(b)$ for any $a,b\ge 0$. The property $f(0)=0$ is immediate from \eqref{eq:fInt}. Monotonicity follows because $\varphi(r) >0 $ and $g(r) \ge 1/2$ by construction.
In particular, we get concavity since for $r\in (0,\infty) \setminus \{R_1,R_2\}$ it holds that
\begin{equation}
   \label{eq:fsecond}
 f''(r) = -h'(r) f'(r) - \frac{\xi}{4}\Ind(r < R_1)-\frac{\beta}{4} \Phi(r) \Ind(r < R_2).
\end{equation}
In particular $f''(r)\le 0$, so that $f'$ is monotonically decreasing.  
The triangle inequality property then follows because
\begin{multline}
\nonumber
f(a+b)=f(a) + \int_{a}^{a+b} f'(s)\d s \\
\le f(a) + \int_{0}^b f'(s)\d s = f(a)+f(b).
\end{multline}

\subsection*{Reflection Coupling}
The key approach from \cite{eberle2019quantitative} is to create an explicit coupling between $\xb_t$ and $\yb_t$ and then bound $\E[\rho(\xb_t,\yb_t)]$. By the definition of the Wasserstein distance from \eqref{eq:wassersteinDef},  $W_\rho(\P_t,\Q_t)\le \E[\rho(\xb_t,\yb_t)]$.
The specific coupling is known as a \emph{reflection coupling} \cite{lindvall1986coupling}.

To define the reflection coupling, let $\taub$ be coupling time:
\begin{equation*}
  \taub = \inf\{t | \xb_t = \yb_t\}.
\end{equation*}
Let $\ub_t = (\xb_t-\yb_t)/\|\xb_t-\yb_t\|$ and let $\ub_t^\star$ be its dual vector so that $\|\ub_t\|^2=\ub_t^\star \ub_t$.  The reflection coupling is given by the following definitions of $\xb_t$ and $\yb_t$:
\begin{align*}
\d\xb_t &= \Hs(\xb_t)\dt + G\d\wb_t - \vb^x_t\d\mub^x(t) \\
\d\yb_t &= \Hs(\yb_t)\dt + (I_n-2\ub_t \ub_t^\star \Ind(t < \taub))G\d\wb_t - \vb_t^y\d\mub_t^y(t).
\end{align*}
Here $-\int_0^t \vb_s^x\d\mub^x(s)$ and $-\int_0^t \vb_s^y\d\mub^y(s)$ are the unique bounded-variation reflection processes that ensure that $\xb_t$ and $\yb_t$ remain in $\Xc$.

Note that $\xb_t=\yb_t$ for all $t\ge \taub$. 

\begin{remark}
Unfortunately, ``reflection'' has two unrelated meanings in this setup: 1) reflecting the processes to remain within $\Xc$ and 2) reflecting the Brownian motion via $I_n-2\ub_t\ub_t^\star$ to couple $\xb_t$ and $\yb_t$.  
\end{remark}

\subsection*{The Supermartingale Property}
Here we show that $e^{at}\rho(\xb_t,\yb_t)$ is a supermartingale.
This is an exercise in non-smooth It\^o calculus for continuous semimartingales. A good reference is Chapter 29 of \cite{kallenberg2021foundations}. 
In the discussion below, $\mb_t$ will denote a local martingale. We use this notation in several places to denote different processes. The specific value of $\mb_t$ of will not matter, since it will have zero mean.

Let $\zb_t = \xb_t-\yb_t$ and $\rb_t = \|\zb_t\|$. We will evaluate $d\rb_t$ via It\^o's formula when $\rb_t > 0$. (Note that $\rb_t >0$ if and only if $t < \taub$.) The method is similar to the proof of Theorem 2.1 in \cite{eberle2019quantitative}, but here extra reflection terms, $- \vb^x_t d\mub^x(t)$ and $- \vb^y_t d\mub^y(t)$,  appear.

For any $z\in \R^n$, let $u = z/\|z\|$. Then we have the Taylor expansion:
$$
\|z+\delta\| = \|z\| + \langle u,\delta \rangle + \frac{1}{2\|z\|}\langle (I_n-uu^\star)\delta,\delta \rangle + o(\|\delta\|^2).
$$
It\^o's formula for continuous-semimartingales gives
\begin{multline}
\nonumber
\d\rb_t = \langle \ub_t, \Hs(\xb_t)-\Hs(\yb_t)\rangle\dt +\langle \ub_t,\vb_t^y \d\mub^y(t)-\vb_t^x \d\mub^x(t)\rangle \\
+\frac{1}{2\|\zb_t\|}\tr\left(
(I_n-\ub_t \ub_t^\star)4 \ub_t \ub_t^\star G G^\star \ub_t \ub_t^\star
\right)\dt + 2\langle \ub_t,G\d\wb_t\rangle.
\end{multline}
  
A significantly simpler upper bound on $\d\rb_t$ can be obtained.
Recall that 
$\ub_t = (\xb_t-\yb_t)/\|\zb_t\|$, $\vb_t^x \in N_{\Xc}(\xb_t)$, and $\vb_t^y \in N_{\Xc}(\yb_t)$. It follows from the definition of the normal cone that $\langle \ub_t ,\vb_t^y\rangle \le 0$ and $-\langle \ub_t ,\vb_t^x\rangle \le 0$. Since $\mub^x$ and $\mub^y$ are non-negative measures, it follows that the corresponding terms can be bounded above by $0$. Furthermore, the trace term vanishes, thus giving:
\begin{equation}
\label{eq:drBound}
d\rb_t \le \langle \ub_t, \Hs(\xb_t)-\Hs(\yb_t)\rangle\dt + \d\mb_t. 
\end{equation}
Recall that $\mb_t$ denotes a zero-mean local martingale.

Now we aim to bound $\d f(\rb_t)$. This part is similar to the analysis from \cite{eberle2019quantitative}. The main distinction is that \cite{eberle2019quantitative} examines the special case of $G=I_n$, while the more general case leads to slightly more complex formulas. We extend $f$ to the entire real line by setting $f(r)=r$ for $r\le 0$. Note that $f$ is a non-smooth concave function which is left-differentiable everywhere.  Denote the left-derivative by $f_{-}'(r)$. Let $\Lb^r_{t}$ denote the right-continuous local time of $\rb_t$. Then Theorem 29.5 of \cite{kallenberg2021foundations} states that
\begin{equation}
\label{eq:genIto}
f(\rb_t)-f(\rb_0)=\int_0^t f_{-}'(\rb_s)\d\rb_s + \frac{1}{2}\int_{-\infty}^{\infty} \Lb^r_t  \d f_{-}'(r),
\end{equation}
where the integral on the right is a Stieltjes integral. (Specifically, it is an integral with respect to the signed measure $\nu$ defined by $\nu([a,b)) = f_{-}'(b)-f_{-}'(a)$ for $a < b$.)

Theorem 29.5 of \cite{kallenberg2021foundations} also states that for any measurable function $b(r)$, we have
\begin{equation}
\label{eq:quadVarToL}
\int_0^t b(\rb_s)\d[\rb]_s = \int_{-\infty}^{\infty} b(r) \Lb_t^r \d r,
\end{equation}
almost surely. 
Here $[\rb]_t = 4\int_0^t\|G^\star \ub_s\|^2 \d s$ is the quadratic variation. Since $G$ is assumed to be invertible, it follows that  $4\sigma_{\min}(G)^2 t \le [\rb]_t \le 4 \|G\|_2^2 t$, where $\sigma_{\min}(G)>0$ is the smallest singular value of $G$. 

Note that $f$ is twice continuously differentiable except at $\{R_1,R_2\}$. For compact notation, set $\zeta=4\sigma_{\min}(G)^2$. Then \eqref{eq:quadVarToL} implies that $\rb_t$  spends zero time at $\{R_1,R_2\}$ almost surely since:
\begin{align}
\nonumber
\MoveEqLeft
\int_0^t \Ind(\rb_s \in \{R_1,R_2\})\d s \le \zeta^{-1} \int_0^t \Ind(\rb_s \in \{R_1,R_2\}) \d[\rb]_s \\
\label{eq:noTime}
&= \zeta^{-1} \int_{-\infty}^{\infty} \Ind(r\in \{R_1,R_2\}) \Lb_t^r\d r = 0.
\end{align}
The last equality follows because $\Lb_t^r$ is right continuous in $r$. More generally, this argument shows that $\rb_t$ spends zero time in any finite collection of points. 

The local time is non-negative, and the measure defined by $\nu([a,b)) = f_{-}'(b)-f_{-}'(a)$ is non-positive. Thus, we get the following inequalities:
\begin{align}
\nonumber
\int_{-\infty}^\infty \Lb_t^r\d f_{-}'(r) &\le   \int_{-\infty}^\infty \Ind(r\notin \{R_1,R_2\})\Lb_t^r\d f_{-}'(r) \\
\nonumber
                                       &= \int_{-\infty}^{\infty} \Ind(r\notin \{R_1,R_2\}) f''(r) \Lb_t^r\d r \\
\nonumber
                                       &=\int_0^t  \Ind(\rb_s\notin \{R_1,R_2\}) f''(\rb_s)\d[\rb]_s \\
\label{eq:localBound}
&\le \zeta\int_0^t f''(\rb_s)\d s.
\end{align}
The last inequality follows from the bounds on the quadratic variation, non-negativity of $f''$ and the fact that $\rb_s$ spends zero time in $\{R_1,R_2\}$ almost surely. 

The preceding argument then implies that for almost all $t\in [0,\taub)$, the following holds almost surely:
\begin{align*}
\MoveEqLeft[0]
\d f(\rb_t)\\
&\le [f'(\rb_t)  \langle \ub_t, \Hs(\xb_t)-\Hs(\yb_t)\rangle + \frac{1}{2}\zeta f''(\rb_t)]\dt + \d\mb_t \\
&\le [f'(\rb_t) (\kappa(r)r + \alpha(\|\xb_t\|+\|\yb_t\|) + \frac{1}{2}\zeta f''(\rb_t)]\dt + \d\mb_t \\
&\le \left[f'(\rb_t)\alpha(\|\xb_t\|+\|\yb_t\|-2M) \right.\\
&\left.-\frac{\xi\zeta}{8}\Ind(r < R_1)-\frac{\beta\zeta}{8}f(r) \Ind(r<R_2)\right]\dt + \d\mb_t.
\end{align*}
The first inequality combined \eqref{eq:drBound}, \eqref{eq:genIto}, and \eqref{eq:localBound}. The second inequality uses the one-sided growth bound \eqref{eq:oneSided}. The third inequality uses \eqref{eq:fsecond}, combined with the following facts, which hold by construction:
\begin{align*}
  f(r) &\le \Phi(r) \\
  h'(r) &= \frac{2}{\zeta}\left(r\kappa(r)+2\alpha M \right).
  \\
  a &\le \beta \zeta / 8 =  \beta \sigma_{\min}(G)/2.
\end{align*}
The fact that $f(r)\le \Phi(r)$ can be deduced from the definitions of $\Phi$ and $f$ from \eqref{eq:fDef}, since $g(r) \in [1/2,1]$. 

Using the Lyapunov assumption, we have
\begin{align*}
\d\Vs(\xb_t)&\le (C-\lambda \Vs(\xb_t))\dt + \d\mb_t \\
\d\Vs(\yb_t)&\le (C-\lambda \Vs(\yb_t))\dt + \d\mb_t.
\end{align*}

Now, we wish to bound $\d(\Vs(\xb_t) \lor \phi(M))$. Recall that $\Vs(x)=\phi(\|x\|)$ and $\phi(r)$ is strictly monotonically increasing.
We will follow a similar strategy as used to bound $\d f(\rb_t)$. 
For compact notation, set $\Vb_t = \Vs(\xb_t)$, let $\hat \Lb_t^r$ be the right-continuous local time of $\Vb_t$,  and let $F(r) = \max\{r,\phi(M)\}$. Note that $F$ is convex and smooth for $r\ne \phi(M)$. Furthermore,  \eqref{eq:noTime} implies that $\rb_t$ spends zero time at $M$, and thus $\Vb_t$ spends zero time at $\phi(M)$. 

Note
that the measure defined by $\hat\nu([a,b)) = F_{-}'(b)-F_{-}'(a)$ is a Dirac delta centered at $\phi(M)$.  Thus, the calculation analogous to \eqref{eq:genIto} gives
$$
F(\Vb_t)-F(\Vb_0)=\int_0^t F_{-}'(\Vb_s)\d\Vb_s + \frac{1}{2} \hat \Lb_t^{\phi(M)}.
$$

Note that the local time, $\Lb_t^{\phi(M)}$ only changes at the times when $\Vb_t = \phi(M)$. See \cite{kallenberg2021foundations}. 
Indeed, the local time is defined by:
\begin{multline}
\nonumber
\hat\Lb^{\phi(M)}_t = |\Vb_t-\phi(M)|-|\Vb_0-\phi(M)| \\- \int_0^t\left(\Ind(\Vb_s>\phi(M))-\Ind(\Vb_s\le \phi(M))\right)\d\Vb_s
\end{multline}
Thus, the local time remains unchanged on intervals in which $\Vb_t > \phi(M)$ or $\Vb_t < \phi(M)$. 
Then, since the amount of time $\Vb_t$ spends at $\phi(M)$ is zero almost surely, we have that for almost all $t< \taub$:
\begin{align*}
\d F(\Vb_t) &= \Ind(\Vb_t > \phi(M)) \d\Vb_t  \\
        &= \Ind(\|\xb_t\| > M) \d\Vs(\xb_t) \\
        &\le \Ind(\|\xb_t\| > M)(C-\lambda \Vs(\xb_t)) + \d\mb_t \\
        &\le -\Ind(\|\xb_t\| > M)(\alpha \|\xb_t\|+\lambda \Vs(\xb_t)/2) + \d\mb_t .
\end{align*}
The first inequality is due to the Lyapunov assumption, while the second inequality follows from the assumption from \eqref{eq:linGrowth}.

Finally, all of the differentials can be combined to give for almost all $t< \taub$:
\begin{multline}
  \nonumber
  d(e^{at} \rho(\xb_t,\yb_t)) \le e^{at}\Big[
    af(\rb_t) -\frac{\beta\zeta}{8}f(\rb_t)\Ind(\rb_t < R_2)+\\
  \left.
    f'(\rb_t)\alpha(\|\xb_t\|+\|\yb_t\|-2M)-\frac{\xi\zeta}{8}\Ind(\rb_t < R_1) \right.\\
  \left.
    +\gamma (2C +(a-\lambda) \Vs(\xb_t) +(a- \lambda) \Vs(\yb_t))\right. \\
  \left.
    +\Ind(\|\xb_t\| > M)((a-\lambda/2)\Vs(\xb_t)-\alpha \|\xb_t\|) \right. \\
   +\Ind(\|\xb_t\| > M)((a-\lambda/2)\Vs(\xb_t)-\alpha \|\xb_t\|)
  \Big]\dt + \d\mb_t.
\end{multline}
Now we will explain why all of the positive terms inside the brackets get canceled by an appropriate negative term.

Note that for $\rb_t < R_2$, the $af(\rb_t)$ term is canceled by  $\frac{\beta\zeta}{8}f(\rb_t)$, since $a\le \beta\eta/8 = \beta \sigma_{\min}(G)/2$. When $\rb_t \ge R_2$, Lemma~2.1 of \cite{eberle2019quantitative} shows that $af(\rb_t)$ is canceled by $\Ac\Vs(\xb_t)+\Ac\Vs(\yb_t)$.

In this case, the triangle inequality implies that $2\max \{\|\xb_t\|,\|\yb_t\|\} \ge \rb_t \ge M/2$. Without loss of generality, assume that $\|\xb_t\|\ge \|\yb_t\|$. Then
$a\le \lambda/2$, combined with \eqref{eq:linGrowth} imply
\begin{align*}
\MoveEqLeft
\gamma(2C+(a-\lambda)\Vs(\xb_t)+(a-\lambda)\Vs(\yb_t)) \\
& \le
 \gamma(2C-(\lambda/2)\Vs(\xb_t)) \\
&\le -\gamma 4C\|\xb_t\| \\
&\le -\frac{\xi \sigma_{\min}(G)^2}{2} \rb_t \\
&\le -a \rb_t
\end{align*}
This cancels $af(\rb_t)$ since $f(\rb_t)\le \rb_t$. 

Now consider the case that $f'(\rb_t)\alpha (\|\xb_t\|+\|\yb_t\|-2M) >0$. Then, since $f'(\rb_t)\ge 0$, at least one of $\|\xb_t\| > M$ or $\|\yb_t\| > M$ must hold. If $\|\xb_t\| > M$, then since $f'(\rb_t)\le 1$, the corresponding term is canceled by
\begin{align*}
  \Ind(\|\xb_t\|>M)((a-\lambda/2)\Vs(\xb_t)-\alpha \|\xb_t\|) \le -\alpha \|\xb_t\|.
\end{align*}
A similar cancellation occurs if $\|\yb_t\| > M$.

The only remaining positive term is now $\gamma 2C = \xi\sigma_{\min}(G)^2/2=\xi\zeta/8 $. In the case that $\rb_t < R_1$, this term is canceled by $-(\xi\zeta/8)\Ind(\rb_t < R_1)$. When $\rb_t\ge R_1$, the definition of $R_1$ along with $a \le \lambda/2$ imply that
\begin{multline}
  \nonumber
  2C + (a-\lambda) \Vs(\xb_t)+(a-\lambda)\Vs(\yb_t)
  \\ \le 2C-(\lambda/2)(\Vs(\xb_t)+\Vs(\yb_t))\le 0.
\end{multline}
Thus, this term is canceled as well.

We have shown that $\d(e^{at}\rho(\xb_t,\yb_t)) \le \d\mb_t$.

\subsection*{Localization}
The last step in proving the theorem requires ruling out the possibility that $\E[e^{at}\rho(\xb_t,\yb_t)]=\infty$. This is performed in \cite{eberle2019quantitative} by a localization argument using stopping times $\taub_k = \inf \{t | \|\rb_t\| \le 1/k \ \textrm{ or } \max\{\|\xb_t\|,\|\yb_t\|\} \ge k\}$ with $\taub_k \to \taub$. The argument works without change in this setting. 
\hfill\QED

\end{document}